\documentclass{amsart}
\usepackage{amsmath, amssymb, latexsym}

\newtheorem{theorem}{Theorem}[section]
\newtheorem{question}[theorem]{Question}
\newtheorem{lemma}[theorem]{Lemma}
\newtheorem{corollary}[theorem]{Corollary}
\newtheorem{definition}[theorem]{Definition}
\newcommand{\mn}{\par\medskip\noindent}
\newcommand{\cB}{\mathcal B}
\newcommand{\Tor}{{\rm Tor\,}}
\def\Cz{\mathbb{C}}
\def\Nz{\mathbb{N}}
\newcommand{\hot}{\mathbin{\overline{\otimes}}}

\title{$L^2$-Betti numbers for subfactors}
\author{Andreas Thom}
\address{Andreas Thom, Mathematisches Institut der Universit\"at G\"ottingen,
Bunsenstr. 3-5, D-37073 G\"ottingen, Germany}
\email{thom@uni-math.gwdg.de}
\urladdr{http://www.uni-math.gwdg.de/thom}

\begin{document}

\begin{abstract}
We study $L^2$-Betti numbers for von Neumann algebras, as defined by D. Shlyakhtenko and A. Connes in \cite{CS}, in the presence 
of a bi-finite correspondence and prove a proportionality formula. 
\end{abstract}

\maketitle

\section{Introduction}

In \cite{CS}, A. Connes and D. Shlyakhtenko define $L^2$-Betti numbers for all complex tracial $*$-algebras $(A,\tau)$ which satisfy a 
certain positivity and boundedness criterion. 
\mn
Let $L^2(A,\tau)$ be the Hilbert space completion of the pre-Hilbert 
space $A$ with inner product $(x,y)_A=\tau(y^*x)$. The boundedness criterion ensures that $A$ acts as bounded multiplication operators on $L^2(A,\tau)$ and so the enveloping von Neumann algebra $M=W^*(A)\subset \cB(L^2(A,\tau))$ exists. Since $M$ carries a faithful trace, it is a finite von Neumann algebra.
The definition of the $k$-th $L^2$-Betti number of $(A,\tau)$ is now as follows:
\[\beta^{(2)}_k(A,\tau)= \dim_{M \hot M^o} \Tor^{A \otimes A^o}_k(M \hot M^o,A) \in [0,\infty].\]

Here, $M\hot M^o$ is naturally seen as right $A \otimes A^o$-module via the inclusion $A \otimes A^o \subset M \hot M^o$. Note that $M \hot M^o$ carries also
a commuting left module structure over $M \hot M^o$. Since the actions commute and $\Tor_k^{A \otimes A^o}(?,A)$ is functorial, $\Tor_k^{A \otimes A^o}(M \hot M^o,A)$ inherits a left-module structure over $M \hot M^o$ with respect to which one can take a dimension. The dimension function to be used is the generalized 
dimension function of W. L\"uck, see \cite{L} on pp. 237.
\mn
The definition is modelled to give a generalization of $L^2$-Betti numbers of discrete groups. Indeed Proposition $2.3$ in \cite{CS} shows that for a
discrete group $\Gamma$ 
\[ \beta^{(2)}_k(\Gamma) = \beta^{(2)}_k(\Cz \Gamma,\tau), \quad \forall k \geq 0,\]
where $\beta^{(2)}_k(\Gamma)$ denotes the $k$-th $L^2$-Betti number of the discrete group $\Gamma$, as studied by J. Cheeger and M. Gromov in \cite{CG}.
\mn
In this article we will study the relation of $L^2$-Betti numbers for $II_1$-factors $M$ and $N$ in the presence of a bi-finite involutive bi-module, relating $M$ and $N$. We
prove a proportionality formula relating the $L^2$-Betti numbers of $M$ and $N$.
\mn
Given a von Neumann algebra $M$, we denote its opposite algebra by $M^o$. $M$-$N$ bi-modules are freely identified with left $M \otimes N^o$-modules. 

\section{Involutive bi-modules}

Let us first set up some notation which we will need in the sequel.

\begin{definition}
Let $M$ and $N$ be von Neumann algebras. Let $L$ be a $M \otimes N^o$-module. The module $L$ is called
bi-finite if it is finitely generated and projective as $M$-module and as $N^o$-module.
\end{definition}

Let $L$ be a bi-finite $M \otimes N^o$-module. As a $M$-module, $L$ is finitely generated and projective and hence isomorphic 
to $M^{\oplus n}p$, for some $n \in \Nz$ and
some projection $p \in M_n(M)$. We have that $$\hom_M(L,L) = (pM_n(M)p)^o$$ and since $N^o$ acts by $M$-module homomorphisms,
we get a unital homomorphism $\phi: N \to pM_n(M)p$. We call $L$ \textit{involutive}, if for some (and hence for any) choice of $n$ and $p$, 
$\phi$ is a $*$-homomorphism.

The requirement of having a $*$-homomorphism is not automatic, however, we know of no interesting examples where this is not 
the case. To resolve this issue completely it would be interesting to know the answer to the following question.

\begin{question}
Let $\phi:N \to M$ be a homomorphism of $II_1$-factors. Let us assume that $M$ is finitely generated and projective as $N$-module. 
Is $\phi$ conjugate to a $*$-homomorphism?
\end{question}

Whereas it is easy to see that such a $\phi$ has to be bounded, we did not succeed in showing that it has to be completely bounded.
It was pointed out to us by D. Shlyakhtenko that, if one could show this, one might be able to use the work of 
G. Pisier on the similarity problem to answer the question affirmatively.
\mn
A source of bi-finite involutive $M\otimes N^o$-modules are $M$-$N$-correspondences. In case of finite 
generation and finite index, their bi-modules of bounded vectors give rise to bi-finite involutive 
$M \otimes N^o$-modules. This will become apparrent in the sequel.
 
\begin{definition}[see \cite{PP}]\label{ppb}
Let $N \subset M$ be a sub-factor of finite index. Denote by $n$ the integer part of $[M:N]$ and let $\alpha = [M:N] -n$. 
Denote by $E: M \to N$ the trace-preserving conditional expectation. 
A Pimsner-Popa basis for the inclusion $N \subset M$ is a finite ordered set $\{m_1, \dots, m_{n+1}\} \subset M$, which
satisfies:
\begin{itemize}
\item $E(m_j m_k^*) = 0, \quad \forall j,k \in \{1, \dots, n+1\}: j \neq k,$
\item $E(m_j m_j^*) = 1, \quad \forall j \in \{1, \dots n\}$
\item $E(m_{n+1} m_{n+1}^*)= p, \quad \mbox{with} \quad p=p^2, \tau(p)= \alpha$
\item $M = \oplus_{i=1}^n N m_i \oplus Npm_{n+1}$, and
\item $\|m_j\| \leq [M:N]^{1/2}, \quad \forall j \in \{1, \dots, n+1\}.$
\end{itemize}
\end{definition}

It is proved in \cite{PP}, that every sub-factor of finite index admits a Pimsner-Popa basis. Note that in particular $M$ is finitely
generated and projective as $N$-module with dimension $[M:N]$. It follows that $M$, seen as a $M \otimes N^o$-module, is bi-finite and involutive.
Conversely, let $N \subset M$ be a sub-factor and assume that $M$ is finitely generated and projective as $N$-module. Then the sub-factor has finite index.
Indeed, this is almost the definition since $\dim_N L^2(M,\tau) = \dim_N M$.

\begin{lemma} \label{basis}
Let $M_0$ be a $II_1$-factor and $N \subset M$ be a sub-factor of finite index.
The canonical map
$$\phi: (M_0 \hot N) \otimes_{M_0 \otimes N} (M_0 \otimes M) \to M_0 \hot M$$ is an isomorphism of left $M_0\hot N$-modules.
\end{lemma}
\begin{proof}
Let $\{m_i, 1 \leq i \leq n\}$ be a Pimsner-Popa basis for the inclusion $N \subset M$. It follows easily that $\{1_{M_0} \otimes m_i, 1 \leq i \leq n\}$
is a Pimsner-Popa basis for the subfactor $M_0 \hot N \subset M_0 \hot M$. The result is now obvious.
\end{proof}

The following two lemmas are easy consequences of the properties of the dimension function. For a proof see \cite{Jones}.

\begin{lemma}\label{comp}
Let $M$ be a $II_1$-factor and let $p$ be a projection in $M_n(M)$ of un-normalized trace $t$. Let $L$ be a module over $M$.
\[\dim_{pM_n(M)p} \, pL^{\oplus n} = \frac1t \,\dim_M \,L.\]
\end{lemma}

\begin{lemma}
Let $N \subset M$ be a subfactor of finite index. For every $M$-module $L$, the following equality holds:
\begin{equation} \label{dim} \dim_N \,L = [M:N] \cdot \dim_M \,L.\end{equation}
\end{lemma}
\section{Main result}

\begin{theorem} \label{main} Let $M$ and $N$ be $II_1$-factors. 
Let $L$ be a bi-finite involutive $M \otimes N^o$-module. Then the following equality holds.
$$\frac{\beta^{(2)}_k(M)}{\dim_M L} = \frac{\beta^{(2)}_k(N)}{\dim_{N^o} L}, \quad \forall k\in \Nz.$$
\end{theorem}
\begin{proof}

As a $M$-module, $L$ is isomorphic to  $M^{\oplus n}p$, for some projection $p \in M_n(M)$. We set $t=n \,\tau(p)$, i.e. $t$ is the un-normalized trace
of $p$. 
Note that in this way, $L$ becomes a $M \otimes (pM_n(M)p)^o$-module.
Since $L$ is involutive, the right $N$-module structure is described by a trace preserving $*$-homomorphism $$\phi: N \to pM_n(M)p = \tilde{M}.$$ 

Using $\phi$, we identify $N$ with a sub-factor
of $\tilde{M}$.
Note that $\tilde{M}$ is finitely generated and projective as a right $N$-module. Indeed, $\tilde{M}$ is a direct
summand in $M_n(M)p$ which is isomorphic to $L^{\oplus n}$. By the remark before Lemma \ref{basis}, 
this implies that the inclusion $N \subset \tilde{M}$ has finite index.

Now, note that $M \otimes N^o \subset M \otimes \tilde{M}^o$ is a flat ring extension, since $M \otimes \tilde{M}^o$ is projective as right $M \otimes N^o$-module, so that there are isomorphisms
\begin{eqnarray} \label{eq0}
\Tor_k^{M \otimes N^o}(M \hot N^o , L) &\cong& \Tor_k^{M \otimes \tilde{M}^o}(M \hot N^o \otimes_{M \otimes N^o} M \otimes \tilde{M}^o, L)\nonumber\\
&\cong& \Tor_k^{M \otimes \tilde{M}^o}(M \hot \tilde{M}^o,L) \quad \mbox{(by \ref{basis})}. \end{eqnarray}
of left $M \hot N^o$-modules. The first isomorphism is the canonical isomorphism of flat base change, see Proposition $3.2.9$ in \cite{W}.
\mn
A computation similar to Theorem $2.4$ in \cite{CS} shows that
\begin{eqnarray}\label{eq1}\dim_{M \hot \tilde{M}^o} \Tor_k^{M \otimes \tilde{M}^o}(M\hot \tilde{M}^o,L) = t^{-1} \beta^{(2)}_k(M).\end{eqnarray}
Indeed, we can define a functor from $M \otimes M^o$-modules to $M \otimes \tilde{M}^o$-modules 
that sends a module $K$ to $ K \otimes_M L = K^{\oplus n}p = (1 \otimes p^o) K^{\oplus n}$. One easily
shows that it is exact and maps finitely generated projective modules to finitely generated projective modules. 
The image of the $M \otimes M^o$-module $M$ is the $M \otimes \tilde{M}^o$-module $L$. 
We conclude that
$$(1\otimes p^o)\left( \Tor_k^{M \otimes M^o}(M \hot M^o,M)\right)^{\oplus n} \cong \Tor_k^{M \otimes \tilde{M}^o}(M \hot \tilde{M}^o,L)$$
Equation \ref{eq1} follows now from Lemma \ref{comp}.
\mn
Since $[M \hot \tilde{M}^o: M \hot N^o] = [\tilde{M}: N],$ we get by equation (\ref{dim}), that
\begin{eqnarray} \label{eq3b} \dim_{M \hot N^o} \Tor_k^{M \otimes N^o}(M \hot N^o,L)  &=& \dim_{M \hot N^o} \Tor_k^{M \otimes \tilde{M}^o}
(M \hot \tilde{M}^o,L) \quad 
\mbox{(by (\ref{eq0}))} \nonumber\\ 
&=& [\tilde{M}:N] \, \dim_{M \hot \tilde{M}^o} \Tor_k^{M \otimes \tilde{M}^o}(M \hot \tilde{M}^o,L) \nonumber \\
&=& t^{-1} \, [\tilde{M}:N] \, \beta^{(2)}_k(M) \quad \mbox{(by (\ref{eq1}))}. \end{eqnarray}

Now, $\dim_{\tilde{M}^o} L = t^{-1}$ and thus
\begin{eqnarray} \label{eq4} \dim_{N^o} L &=& [\tilde{M}:N] \, \dim_{\tilde{M}^o} L \quad \mbox{(by (\ref{dim}))}  \nonumber\\
&=& t^{-1} \, [\tilde{M}: N] \end{eqnarray}
Equations (\ref{eq3b}) and (\ref{eq4}) imply
\begin{equation} \label{done} \dim_{M \hot N^o} \Tor_k^{M \otimes N^o}(M \hot N^o,L) = \dim_{N^o} L \cdot \beta^{(2)}_k(M). \end{equation}
The claim follows by symmetry of the left hand side of equation (\ref{done}).
\end{proof}

\begin{corollary}
Let $N \subset M$ be a subfactor of finite index. The following relation holds.

$$ [M:N] \cdot \beta^{(2)}_k(M) = \beta^{(2)}_k(N), \quad \forall k \geq 0$$
\end{corollary}
\begin{proof} The $N \otimes M^o$-module $M$ is bi-finite by Lemma \ref{ppb} and involutive. Therefore the corollary is
implied by the theorem.
\end{proof}

There is an independent proof by L\"uck of the above corollary in \cite{L2} in the case of inclusions $L H \subset L G$, coming from 
an inclusion of a normal subgroup $H \subset G$ of finite index.

\end{document}